\documentclass[11pt]{article}
\usepackage{amssymb,amsmath,latexsym}
\def\msy{\mathbb}
\newcommand{\qed}{\unskip\nobreak\hfill\mbox{ $\Box$}\par}

\def\bbbc{{\msy C}}
\def\bbbq{{\msy Q}}

\def\bbbz{{\msy Z}}
\def\bbba{\overline{\bbbq}}

\def\fc#1#2{\left[{#1\atop #2}\right]}

\def\gcd{{\rm gcd}}
\def\lcm{{\rm lcm}}

\def\is{\equiv}
\def\mod#1{{(\rm mod\ }#1)}
\newtheorem{theorem}[subsection]{Theorem}

\newtheorem{lemma}[subsection]{Lemma}

\newtheorem{corollary}[subsection]{Corollary}

\newtheorem{proposition}[subsection]{Proposition}

\begin{document}

\title{Irrationality of some $p$-adic $L$-values}
\author{Frits Beukers\\
Department of Mathematics, University of Utrecht, P.O.\ Box 80.010,
3508 TA Utrecht, The~Netherlands}
\date{June 16, 2006}
\maketitle

\begin{abstract}
We give a proof of the irrationality of the $p$-adic
zeta-values $\zeta_p(k)$ for $p=2,3$ and $k=2,3$. Such
results were recently obtained by F.Calegari as an application
of overconvergent $p$-adic modular forms. In this paper we
present an approach using classical continued fractions
discovered by Stieltjes. In addition we show irrationality of some other
$p$-adic $L$-series values, and values of the $p$-adic
Hurwitz zeta-function. 
\end{abstract}

\section{Introduction}
The arithmetic nature of values of Dirichlet $L$-series at
integer points $r>1$ is still a subject with many unanswered
questions. It is
classically known that if a Dirichlet character $\chi:\bbbz\to 
\bbbc$ has the same parity as $r$, the number $L(r,\chi)$ is
an algebraic multiple of $\pi^r$, hence transcendental. When
$\chi$ has parity opposite from $r$, the matter is quite different.
The only such value known to be irrational is $\zeta(3)$ as
R.Ap\'ery first proved in 1978. From later work by Rivoal and Ball
\cite{RivBal01} it follows that $\zeta(2n+1)$ is irrational for infinitely many
$n$, and W.Zudilin \cite{Zud04} showed recently that at least one among
$\zeta(5),\zeta(7),\zeta(9),\zeta(11)$ is irrational. We also recall
analogous statements for $L$-values with the odd character modulo 4 in
\cite{RivZud03}.

Although there have been many attempts to generalise Ap\'ery's
original irrationality proof to higher zeta-values, all have failed
due to the absence of convenient miracles which did occur in the
case of $\zeta(3)$. One such attempt was made by the present author
\cite{Beu87} through the use of elementary modular forms. Although the approach
looked elegant it provided no new significant results.
Ever since then the method
has lain dormant with no new applications.

In a recent, very remarkable and beautiful paper, Frank Calegari \cite{Cal05}
managed to establish further irrationality results using modular
forms. However, the numbers involved are values of Leopoldt-Kubota
$p$-adic $L$-series. 
For example, Calegari managed to prove irrationality
of $2$-adic $\zeta(2)$ and $2$- and $3$-adic $\zeta(3)$. The underlying
mechanism is the overconvergence of certain $p$-adic modular forms,
a subject which has recently attracted renewed attention in connection
with deformation theory of Galois representations. 

Since overconvergent modular form theory is an advanced subject
I tried to reverse engineer the results of Calegari in order to toss
out the use of modular forms and find a more classical approach.
This turns out to be possible. We show irrationality of a large family
of $p$-adic numbers, some of which turn out to be values of
$p$-adic $L$-series at the points 2 or 3.
In Theorems \ref{mainI}, \ref{mainII}
and \ref{mainIII} one finds
the main results of this paper. Incidently we note that in Calegari's
paper irrationality of the $2$-adic Catalan constant is mentioned. 
Although the term 'Catalan constant' is perfectly reasonable,
it does not correspond with the Kubota-Leopoldt $L_2(2,\chi_4)$
where $\chi_4$ is the odd Dirichlet character modulo $4$.
It is well-known that Kubota-Leopoldt $L$-series with odd character
vanish identically. The work of Calegari actually entails irrationality of the
Kubota-Leopoldt $\zeta_2(2)$. The difference
is due to the extra Teichm\"uller character which occurs in the
definition of the Kubota-Leopoldt $L$-functions.
In \cite{Cal05} the irrationality of $2$- and $3$-adic $\zeta(3)$ is also shown.

In Sections \ref{padeI}, \ref{padeII} and \ref{padeIII} we shall discuss
$p$-adic irrationality results proved using Pad\'e approximations to the infinite
Laurent series
\begin{eqnarray*}
\Theta(x)&=&\sum_{n\ge0}t_n(-1/x)^{n+1} \\
R(x)&=&\sum_{n\ge0}B_n(-1/x)^{n+1}\\
T(x)&=&\sum_{n\ge0}(n+1)B_n(-1/x)^{n+2}
\end{eqnarray*}
where the $B_n$ are the Bernoulli numbers and $t_n=(2^{n+1}-2)B_n$.
Continued fraction expansions to $R(x)$ and $T(x)$ were already known
to T.J.Stieltjes in 1890. In a first version of this paper I worked
out the corresponding Pad\'e approximations using hypergeometric
functions, which can be found in this paper. However, it was pointed out to me
by T.Rivoal that the Pad\'e approximations for $R(x)$ and $T(x)$
were also described in a different way in a paper by M.Pr\'evost \cite{Pre96}.
In this paper the author gives an alternative irrationality proof of
Ap\'ery's result $\zeta(2),\zeta(3)\not\in\bbbq$. In a paper by Rivoal
\cite{Riv06} the author makes a similar attempt at proving irrationality of Catalan's
constant. The Pad\'e approximations involved in there are precisely the 
approximations to $\Theta(x)$! Ironically in both \cite{Pre96} and 
\cite{Riv06} the implications
for proving $p$-adic irrationality results are not noted. 

The irrationality results of Calegari are contained in the irrationality
results that were found in the Pad\'e approximation approach sketched above.

We collected the definition and basic properties of $p$-adic $L$-series
in the Appendix of this paper. Throughout we use the conventions
made in Washington's book \cite[Ch.5]{Was97} on cyclotomic fields.

{\bf Acknowledgements}.
I am deeply grateful to Henri Cohen for having me
provided with a proof of Proposition \ref{cohen}, which was a crucial step
in writing up this paper. 
Details of his proof will occur as
exercises in Cohen's forthcoming book on Number Theory.
The proof presented
here is a shorter but less transparent one, derived from Cohen's 
observations. 

Thanks are also due to the authors of the number theory package PARI which
enabled me to numerically verify instances of $p$-adic identities. I also
thank the authors P.Paule, M.Schorn, A.Riese of the Fast Zeilberger package
for Mathematica. Their implementation of Gosper and Zeilberger
summation turned out
to be extremely useful. 

Finally I like to thank T.Rivoal for pointing out the connections with
existing irrationality results.

\section{Arithmetic considerations}
The principle of proving irrationality of a $p$-adic number $\alpha$
is to construct a sequence of rational approximations $p_n/q_n$
which converges $p$-adically to $\alpha$ sufficiently fast. To be
more precise,

\begin{proposition}\label{irrationality}
Let $\alpha$ be a $p$-adic number and
let $p_n,q_n,\ n=0,1,2,\ldots$ be two sequences of integers
such that 
$$\lim_{n\to\infty}\max(|p_n|,|q_n|)|p_n-\alpha q_n|_p=0$$
and $p_n-\alpha q_n\ne0$ infinitely often. Then $\alpha$ is 
irrational.
\end{proposition}

{\bf Proof}. Suppose $\alpha$ is rational, say $A/B$ with $A,B\in
\bbbz$ and $B>0$. Whenever $p_n-(A/B)q_n$ is non-zero we have
trivially, $\max(|p_n|,|q_n|)|p_n-(A/B)q_n|_p\ge 1/\max(|A|,B)$.
Hence the limit as $n\to\infty$ cannot be zero. Thus we conclude
that $\alpha$ is irrational.
\qed
\medskip

We will also need some arithmetic statements about hypergeometric
coefficients.

\begin{lemma}\label{pochhammerrational}
Let $\beta$ be a rational number with the integer
$F\in\bbbz_{>1}$ as denominator.
Then $(\beta)_n/n!$ is a rational number whose denominator
divides $\mu_n(F)$, where
$$\mu_n(F)=F^n \prod_{q|F} q^{[n/(q-1)]}$$
where the product is over all primes $q$ dividing $F$.
Moreover, the number of primes $p$ in the denominator of $(\beta)_n/n!$
is at least $n(r+1/(p-1))-\log n/\log p-1$, where $r$ is defined by 
the relation $|F|_p=p^{-r}$.
\end{lemma}

{\bf Proof}. Let us write $\beta=b/F$ with $b\in\bbbz$. Then
$${(\beta)_n\over n!}={\prod_{k=0}^{n-1}(b+Fk)\over F^n n!}.$$
Let $q$ be a prime. Suppose $q$ divides $F$. Then $q$ does not
divide the product $\prod_k(b+Fk)$ and the number of primes $q$
in the denominator is the number of primes $q$ in $F^n n!$.
The number of primes $q$ in $n!$ equals
$$[n/q]+[n/q^2]+[n/q^3]+\cdots$$
which is bounded above by
$$n/q+n/q^2+n/q^3+\cdots=n/(q-1).$$
This explains the factor $q^{[n/(q-1)]}$ in our assertion.
Morever, we also have the lower bound
$$\sum_{k=1}^{[\log n/\log q]}[n/q^k]\ge
\sum_{k=1}^{[\log n/\log q]}(n/q^k-1)\ge {n-1\over q-1}-{\log n
\over \log q}.$$
This lower bound accounts for the second assertion when $q=p$.

To finish the proof of the first assertion we must show that
$(n!)^{-1}\prod_{k=0}^{n-1}(b+Fk)$ is $q$-adically integral if
$q$ does not divide $F$. This follows easily from the fact that
the number of $0\le k\le n-1$ for which $b+Fk$ is divisible by a
power $q^s$ is always larger or equal than the number of $1\le k\le n$
for which $k$ is divisble by $q^s$.
\qed

\section{Differential equations}

In the next sections we shall consider solutions of linear
differential equations of orders 2 and 3. Here we derive
some generalities on the arithmetic of the coefficients
of the solutions in Taylor series. 

Let $R$ be a domain of characteristic zero with quotient field $Q(R)$.
Consider a differential operator $L_2$ defined by
$$L_2(y):=zp(z)y''+q(z)y'+r(z)y$$
where $p(z),q(z),r(z)\in R[z]$, $p(0)=1$. Suppose there exists
$W_0\in R[[z]]$ such that $W_0(0)=1$ and the logarithmic derivative
of $W_0/z$ equals $-q(z)/zp(z)$. We call $W_0/z$ the Wronskian determinant
of $L_2$. Suppose in addition that the equation
$L_2(y)=0$
has a formal power series $y_0\in R[[z]]$ with $y_0(0)=1$
as solution. Such a solution is determined uniquely since the space
of solutions in $Q(R)[[z]]$ has dimension one.

The operator $L_2$ has a symmetric square $L_3$ which we write as 
$$L_3(y):=z^2P(z)y'''+Q(z)y''+R(z)y'+S(z)y$$
with $P,Q,R,S\in R[z]$, $P(0)=1$. This symmetric square is
characterised by the property that the solution space of $L_3(y)$
is spanned by the squares of the solutions of $L_2(y)=0$.
The equation $L_3(y)=0$ has a unique formal power series solution
with constant term $1$, which is $y_0^2$. 

\begin{proposition}\label{almostintegral}
Let notations and assumptions be as above. Then the inhomogeneous equation 
$L_2(y)=1$ has a unique solution $y_{\rm inhom,2}\in Q(R)[[z]]$
starting with $z+O(z^2)$. Moreover, the $n$-th coefficient
of $y_{\rm inhom,2}$ has denominator dividing $\lcm(1,2,\ldots,n)^2$.

The inhomogeneous equation $L_3(y)=1$ has a unique solution
$y_{\rm inhom,3}\in Q(R)[[z]]$ starting with $z+O(z^2)$. Moreover, the
$n$-th coefficient of $y_{\rm inhom,3}$ has denominator dividing
$\lcm(1,2\ldots,n)^3$.
\end{proposition}

{\bf Proof}.
In this proof we shall use the following identities, which hold
for any $f\in Q(R)[[z]]$,
$$\int_0^z f \log z\ dz=\log z\int_0^z f dz -\int_0^z{1\over z}\int_0^z fdz.$$
and
$$\int_0^z f (\log z)^2\ dz=(\log z)^2\int_0^z f dz -2\log z\int_0^z{1\over z}
\int_0^z fdz + 2\int_0^z{1\over z}\int_0^z{1\over z}\int_0^z fdz.$$
These identities can be shown by (repeated) partial integration.

One easily verifies that a second, independent solution of $L_2(y)=0$
is given by $y_1=y_0\int (W_0/zy_0^2)dz$. The quotient $W_0/zy_0^2$ equals
$1/z$ plus a Taylor series in $R[[z]]$. Integration and multiplication
by $y_0$ then shows that
$y_1=y_0\log z+\tilde{y}_0$ where $\tilde{y}_0\in Q(R)[[z]]$
and whose $n$-th coefficient has denominator dividing
$\lcm(1,2,\ldots,n)$. We choose the constant of integration in such a way
that $\tilde{y}_0(0)=0$.

Note by the way that $y_1'y_0-y_1y_0'=W_0/z$ which is precisely how the
Wronskian should be defined. A straightforward verification shows that
$$y_1\int_0^z {y_0\over pW_0}\ dz-y_0\int_0^z {y_1\over pW_0}\ dz$$
is solution of the inhomogeneous equation $L_2(y)=1$. Now
substitute $y_1=y_0\log z+\tilde{y}_0$. We obtain, using the identity
for $\int_0^z f \log z\ dz$,
$$\tilde{y}_0\int_0^z{y_0\over pW_0}\ dz-y_0\int_0^z 
{\tilde{y}_0\over pW_0}\ dz
+y_0\int_0^z{1\over z}\int_0^z {y_0\over pW_0}\ dz\ dz.$$
We have thus obtained a power series solution of $L_2(y)=1$ and
the assertion about the denominators of the coefficients readily follows.

Another straightforward calculation shows that
$$y_0^2\int_0^z{y_1^2\over W_0^2P}\ dz-2y_0y_1\int_0^z{y_0y_1\over
W_0^2P}dz +y_1^2\int_0^z{y_0^2\over W_0^2P}\ dz$$
is a solution of $L_3(y)=2$. Continuation of our straightforward calculation
using $y_1=y_0\log z+\tilde{y}_0$ shows that this solution equals
$$2y_0^2\int_0^z{1\over z}\int_0^zY\ dz\ dz+2y_0\tilde{y}_0\int_0^zY\ dz
+\tilde{y}_0^2\int_0^z{y_0^2\over W_0^2P}\ dz$$
where
$$Y=-{y_0\tilde{y}_0\over W_0^2P}+
{1\over z}\int_0^z {y_0^2\over W_0^2P}\ dz.$$
The last statement of our Proposition follows in a straightforward
manner. 

\qed

\section{Some identities}
Consider the field of rational functions $\bbbq(x)$ with a discrete valuation
such that $|x|>1$. Denote its completion with respect to that valuation by $K$.
We see that $K$ is the field of formal Laurent series in $1/x$. Our considerations
will take place within this field.
Later we shall substitute $x=a/F$
where $a/F$ is a rational number with $|a/F|_p>1$ and perform an evaluation in
$\bbbq_p$. 

Define, following J.Diamond in \cite{Dia79},
$$\fc{n}{x}={n!\over x(x+1)\cdots(x+n)}.$$

\begin{proposition}\label{basicidentities}
Let $\Theta(x),R(x),T(x)\in K$ be the Taylor series in $1/x$ which we
defined in the introduction.
Then we have the following identities in $K$,
$$\Theta(x)=-\sum_{n=0}^{\infty}\fc{n}{x}\fc{n}{1-x},$$
$$R(x)=-\sum_{k=0}^{\infty}{1\over k+1}\fc{k}{x},$$
$$T(x)=-\sum_{k=0}^{\infty}{1\over k+1}\fc{k}{x}\fc{k}{1-x}.$$
\end{proposition}

It will be the purpose of this section to prove these equalities.
Let us first record the following relations between $\Theta(x),R(x)$ and
$T(x)$ which follow directly from their definition. Namely
$$T(x)=R'(x)\qquad \Theta(x)=R(x/2)-2R(x).$$

For any $A(x)\in K$ we can also consider $A(x+\lambda)$ for any $\lambda\in
\bbbq$ as element of $K$ if we expand $1/(x+\lambda)$ formally in a
power series in $1/x$ again. We use the following important observation.

\begin{lemma}\label{invariance}
Suppose $A(x)\in K$ and suppose there exists a non-zero $\lambda\in\bbbq$ such
that $A(x+\lambda)=A(x)$. Then $A(x)$ is a constant.
\end{lemma}

{\bf Proof}. The equality $A(x+\lambda)=A(x)$ remains true if we subtract
the constant coefficient $a_0$ from $A$. Let us now assume $A(x)-a_0$
that $A(x)-a_0$ is not identically zero.
Then there exists a non-zero integer $n$ and non-zero $a_n$
such that $A(x)-a_0=a_n(1/x)^n+$ higher order terms in $1/x$.
It is straightforward to verify that $A(x+\lambda)-A(x)=-n\lambda a_n
(1/x)^{n+1}+$ higher order terms in $1/x$. This contradicts
$A(x+\lambda)=A(x)$. Hence $A(x)-a_0$ is identically zero. 

\qed
\medskip

We require the following property of Bernoulli-numbers.

\begin{lemma}\label{bernoullishift}
For any $n\ne1$ we have 
$$\sum_{k=0}^nB_k{n\choose k}=B_n.$$
When $n=1$ we have $B_0+B_1=1+B_1$.
\end{lemma}

{\bf Proof}. Recall the definition
$${t\over e^t-1}=\sum_{n\ge0}B_nt^n/n!.$$
Multiplication by $e^t$ gives 
$$t+{t\over e^t-1}=\sum_{n\ge0}\left(\sum_{k=0}^n B_k{n\choose k}\right)t^n/n!$$
Our Lemma follows by comparison of coefficients of $t^n$.

\qed
\medskip

We are now ready to prove the following functional equations.
\begin{proposition}\label{functionalR}
We have the identities
\begin{enumerate}
\item $R(x+1)-R(x)=1/x^2$
\item $R(x)+R(1-x)=0$
\item $R(x)+R(x+1/2)=4R(2x)$
\end{enumerate}
\end{proposition}

{\bf Proof}. 
The first statement follows from 
\begin{eqnarray*}
R(x+1)&=&\sum_{k\ge0}B_k(-1/(x+1))^{k+1}\\
&=&\sum_{k\ge0}B_k\sum_{n=0}^{\infty}{n\choose k}(-1/x)^{n+1}
\end{eqnarray*}
Now interchange the summations to get
\begin{eqnarray*}
R(x+1)&=&\sum_{n\ge0}\sum_{k=0}^n B_k{n\choose k}(-1/x)^{n+1}\\
&=&1/x^2+\sum_{n\ge0}B_n(-1/x)^{n+1}
\end{eqnarray*}
where the last equality follows from Lemma \ref{bernoullishift}.

To show the second statement we use the identity $R(-x)+R(x)=1/x^2$
which follows from the fact that the only odd index $n$ for which
$B_n\ne 0$ is $n=1$. Combining this with the first statement yields
the second statement.

To show the third statement we use Lemma \ref{invariance}. 
Write $A(x)=4R(2x)-R(x)-R(x+1/2)$. Notice that $A(x)$ has constant
term zero and from our first two results we deduce
\begin{eqnarray*}
A(x+1/2)-A(x)&=&4R(2x+1)-4R(2x)-R(x)+R(x+1)\\
&=&4/(4x^2)-1/x^2=0.
\end{eqnarray*}
Hence our Lemma implies that $A(x)$ is identically zero.

\qed
\medskip

For $T(x),\Theta(x)$ there are a few immediate corollaries.
\begin{corollary}\label{identities}
We have
\begin{enumerate}
\item $T(x+1)-T(x)=-2/x^3$
\item $T(x)=T(1-x)$
\item $\Theta(x+1)+\Theta(x)=-2/x^2$
\item $\Theta(x)=R(x/2)-R(x/2+1/2).$
\end{enumerate}
\end{corollary}

{\bf Proof}. The first two statement follow from the first two statements
of Proposition \ref{functionalR} because $T(x)=R'(x)$.

For the third statement we use $\Theta(x)=R(x/2)-2R(x)$ and the
third statement of Proposition \ref{functionalR}. We get
\begin{eqnarray*}
\Theta(x+1)+\Theta(x)&=&R(x/2)+R(x/2+1/2)-2R(x)-2R(x+1)\\
&=&4R(x)-2R(x)-2R(x+1)=-2/x^2.
\end{eqnarray*}

The last statement follows from $\Theta(x)=R(x/2)-2R(x)$ and 
$4R(x)=R(x/2)+R(x/2+1/2)$.

\qed
\medskip

We are now ready to prove Proposition \ref{basicidentities}.
To prove the first identity we set
$$S(x)=\sum\fc{n}{x}\fc{n}{1-x}$$
and show that it satisfies $S(x+1)+S(x)=-2/x^2$. Our assertion then
follows from $S(x+1)-\Theta(x+1)+S(x)-\Theta(x)=0$, hence 
$S(x)-\Theta(x)$ is periodic
with period 2. Application of Lemma \ref{invariance} then shows that
$S(x)-\Theta(x)$ is identically zero. 

By straightforward calculation we find
$$\fc{n}{x}\fc{n}{1-x}+\fc{n}{x+1}\fc{n}{-x}=2\fc{n}{x+1}\fc{n}{1-x}.$$
Using Gosper summation we get
$$\fc{n}{x+1}\fc{n}{1-x}=\Delta_n\left({(n+1-x)(n+1+x)\over x^2}
\fc{n}{1-x}\fc{n}{1+x}\right)$$
where $\Delta_n$ is the forward difference operator $\Delta_n(g)(n)
=g(n+1)-g(n)$.
Now carry out the summation and use telescoping of series to
find that $S(x+1)+S(x)=-2/x^2$.

To prove the second assertion of Proposition \ref{basicidentities}
we denote the summation 
on the right again by $S(x)$. Observe that
$$\fc{k}{x+1}-\fc{k}{x}=-{k+1\over x}\fc{k}{x+1}.$$
Hence 
$$S(x+1)-S(x)=\sum_{k=0}^{\infty}{1\over x}\fc{k}{x+1}.$$
Using Gosper summation one quickly finds that
$$\fc{k}{x+1}=-{1\over x}\Delta_k\left((1+k+x)\fc{k}{x+1}\right).$$
Summation over $k$ then yields
$$S(x+1)-S(x)={1\over x^2}.$$
Hence $R(x)-S(x)$ is periodic with period $1$ and thus identically $0$
according to Lemma \ref{invariance}.

To prove the third assertion of Proposition \ref{basicidentities}
we again denote the righthand side
by $S(x)$. Observe that
$$\fc{k}{x+1}\fc{k}{-x}-\fc{k}{x}\fc{k}{1-x}=
-2{k+1\over x}\fc{k}{1+x}\fc{k}{1-x}.$$
Hence
$$S(x+1)-S(x)=\sum_{k=0}^{\infty}{-2\over x}\fc{k}{1+x}\fc{k}{1-x}.$$
Using Gosper summation one easily finds that
$$\fc{k}{1+x}\fc{k}{1-x}={-1\over x^2}\Delta_k\left((x^2-(k+1)^2)
\fc{k}{x+1}\fc{k}{1-x}\right).$$
Summation over $k$ then yields
$$S(x+1)-S(x)=-{2\over x^3}.$$
Hence $T(x)-S(x)$ is periodic with period $1$ and thus identically $0$
according to Lemma \ref{invariance}.

\qed
\medskip

\section{Some $p$-adic identities}
In the following results we relate $p$-adic values of $R(x),T(x),\Theta(x)$
with some $p$-adic
$L$-series. Let $a/F$ be a rational number whose denominator
$F$ is divisible by $p$. The series obtained from $\Theta(x),R(x),T(x)$ by the
substitution $x=a/F$ are $p$-adically convergent. We denote the $p$-adic
values of these series by $\Theta_p(a/F)$, $R_p(a/F)$ and $T_p(a/F)$.
First of all, it follows in a straightforward manner from the Appendix that
$$R_p(a/F)=-F^2\omega(a)^{-1}H_p(2,a,F)$$
where $H_p$ is the $p$-adic Hurwitz zeta-function and $\omega$ the
Teichm\"uller character modulo $p$. As a Corollary we get expressions
for $\Theta_p(a/F)$ in terms of $p$-adic Hurwitz zeta-function values.

As application we now have,
\begin{proposition}\label{cohen}
Let $\chi_d$ be the primitive even Dirichlet character
modulo $d$. Then,
\begin{enumerate}
\item $\Theta_2(1/2)=-8\zeta_2(2)$
\item $\Theta_2(1/6)=-40\zeta_2(2)$
\item $\Theta_2(1/4)=-16L_2(2,\chi_8)$
\item $\Theta_3(1/3)=-27\zeta_3(2)/2$
\item $\Theta_3(1/6)=-36L_3(2,\chi_{12})$
\end{enumerate}
\end{proposition}

We show how to prove the first assertion. Using Corollary
\ref{identities} (4) we get
$$\Theta_2(1/2)=(R_2(1/4)-R_2(3/4))/2.$$
Using the relation between $R_2(a/F)$ and $H_2(2,a,F)$ sketched above we find
$$\Theta_2(1/2)=-8(H_2(2,1,4)+H_2(2,3,4))=8\zeta_2(2),$$
where the last equality follows from the last formula in the Appendix.

Similarly we can find $p$-adic values of $T(x)$ as $p$-adic
zeta-values. We use the fact that
$$T_p(a/F)=2F^3\omega(a)^{-2}H_p(3,a,F).$$
As a consequence we get the following evaluations.

\begin{proposition}\label{valuesat3}
Let $\chi_d$ be the even primitive Dirichlet character modulo $d$. Then
\begin{enumerate}
\item $T_2(1/4)=4^3\zeta_2(3)$
\item $T_3(1/3)=3^3\zeta_3(3)$
\item $T_5(1/5)=(5^3/2)(\zeta_5(3)-L_5(3,\chi_5))$
\item $T_5(2/5)=(5^3/2)(\zeta_5(3)+L_5(3,\chi_5))$
\item $T_2(1/8)=2^8(\zeta_2(3)-L_2(3,\chi_8))$
\item $T_2(3/8)=2^8(\zeta_2(3)+L_2(3,\chi_8))$
\end{enumerate}
\end{proposition}

As illustration we prove the first equality. First use $T(x)=T(1-x)$
to get $T_2(1/4)=T_2(3/4)$. Then observe,
$$T_2(1/4)=(T_2(1/4)+T_2(3/4))/2=4^3(H_2(3,1,4)+H_2(3,3,4))=4^3\zeta_2(3).$$

\section{Pad\'e approximations I}\label{padeI}
In this section we will prove that $\Theta(x)$ has the follwoing
continued fraction expansion,
$$\Theta(x)=\cfrac{1}{x^2-x+a_1-\cfrac{b_1}{x^2-x+a_2-\cfrac{b_2}{\ddots}}}$$
where $a_n=2n^2-2n+1,b_n=n^4$.
We shall study its convergents and use these to derive irrationality
results. For example, if we substitute $x=1/2$ we obtain a continued
fraction expansion which converges $2$-adically very fast to the $2$-adic
evaluation $\Theta_2(1/2)=-8\zeta_2(2)$. In fact, the convergents of this
continued frcation for $x=1/2$ are precisely Calegari's approximations found
in \cite{Cal05}. 

From the theory of continued fractions it follwos that the convergents
are of the form $V_n/U_n\ n=0,1,2,\ldots$ where $V_n,U_n$ are polynomials of degrees $2n-2, 2n$ respectively. Moreover $U_n,V_n$ satisfy the recurrence
relation
$$U_{n+1}=(2n^2+2n+1-x+x^2)U_n-n^4U_{n-1}.$$
Now substitute $U_n=(n!)^2u_n$ and we get a new recurrence relation
\begin{equation}\label{recurrenceI}
(n+1)^2u_{n+1}=(2n(n+1)+1-x+x^2)u_n-n^2u_{n-1}
\end{equation}
Consider the solutions $q_n(x)$ and $p_n(x)$ given by
\begin{eqnarray*}
q_0(x) &=& 1\\
q_1(x) &=& x^2-x+1\\
q_2(x) &=& (x^4-2x^3+7x^2-6x+4)/4\\
q_3(x) &=& (x^6-3x^5+22x^4-39x^3+85x^2-66x+36)/36\\
&&\cdots
\end{eqnarray*}
and
\begin{eqnarray*}
p_0(x) &=& 0\\
p_1(x) &=& 1\\
p_2(x) &=& (x^2-x+5)/4\\
p_3(x) &=& (x^4-2x^3+19x^2-18x+49)/36\\
&&\cdots
\end{eqnarray*}
Then the sequence of rational functions $p_n(x)/q_n(x)$ are the convergents
of our continued fraction. To determine the $q_n(x)$ we consider the 
generating function
$$y_0(z)=\sum_{n=0}^{\infty}q_n(x)z^n.$$
Due to the recursion relation it is straightforward to see that $y_0(z)$
is a power series solution of the second order linear differential
equation
$$L_2(y)=z(z-1)^2y''+(3z-1)(z-1)y'+(z-1+x(1-x))y=0$$
where the $'$ denotes differentation with respect to $z$.
Power series solutions in $z$ are uniquely determined up to a
scalar factor. Since it is also straightforward to see that
$(1-z)^{x-1}\ _2F_1(x,x,1,z)$ is another such solution we conclude
that
$$y_0(z)=(1-z)^{x-1}\ _2F_1(x,x,1,z).$$
Comparison of coefficients gives us
$$q_n(x)=\sum_{k=0}^n{(1-x)_{n-k}(x)_k^2\over (n-k)!(k!)^2}.$$
Let us also consider the generating function for the $p_n(x)$,
$$y_1(z)=\sum_{n=0}^{\infty}p_n(x)z^n.$$
A straightforward calculation using the recurrence shows that 
$L_2(y_1)=1$.

The problem is now to show that the rational functions $p_n(x)/q_n(x)$
approximate $\Theta(x)$ in $K$. To that end we define for each $n$,
$$\Theta(n,x)=(-1)^n\sum_{k=0}^{\infty}{k\choose n}\fc{k}{x}\fc{k}{1-x}.$$
Notice that $\Theta(0,x)=-\Theta(x)$ via Proposition \ref{basicidentities}.

\begin{proposition}\label{remainderI}
Letting notations be as above, we have for each $n$,
$$p_n(x)-q_n(x)\Theta(x)=\Theta(n,x)=O(1/x^{2n+2})$$
as Laurent series in $1/x$.
\end{proposition}
{\bf Proof}. Letting
$$F(k,n)=(-1)^n{k\choose n}\fc{k}{x}\fc{k}{1-x}$$
the Zeilberger algorithm shows that
$$n^2F(k,n-1)-(-x+x^2+2n^2+2n+1)F(k,n)+(n+1)^2F(k,n+1)=
\Delta_k(F(k,n)(x+k)(k+1-x)).$$
When $n\ge1$ summation over $k$ yields
$$n^2\Theta(n-1,x)-(-x+x^2+2n^2+2n+1)\Theta(n,x)+(n+1)^2\Theta(n+1,x)=0.$$
When $n=0$ we get
$$-(-x+x^2+1)\Theta(0,x)+\Theta(1,x)=1.$$
From this, and the fact that $\Theta(0,x)=-\Theta(x)$ we conclude that
$$p_n(x)-q_n(x)\Theta(x)=\Theta(n,x)$$
for all $n\ge0$.

\qed
\medskip

It was remarked to me by T.Rivoal that these approximations can also
be found in \cite{Riv06} as $P_{2n}(z)$. When we replace the $z$ there by $1-2x$
we obtain the alternative expression,
$$q_n(x)=\sum_{k=0}^n{n\choose k}{-x\choose k}{k-x\choose k},$$
Notice that by taking $x=-n$ we recover Ap\'ery's numbers for
the irrationality of $\zeta(2)$ again. By taking $x=-n+1/2$ one
obtains numbers which play a role in approximations of Catalan's constant
(see \cite{Riv06}).

From \cite{Riv06} we find an explicit formula for $p_n$ (known as $Q_{2n}$ in
\cite{Riv06}),
$$p_n(x)=\sum_{k=1}^n{n\choose k}\sum_{j=1}^k{k-x\choose k-j}
{-x-j\choose k-j}{(-1)^{j-1}\over j^2{k\choose j}^2}.$$

\section{Application I}
\begin{proposition}\label{estimatesI}
Let $p_n,q_n$ be as in the previous section and let $\mu_F(n)$ be as in Lemma 
\ref{pochhammerrational}.
Then, 
\begin{enumerate}
\item For every $n$ the number $q_n(a/F)$ is rational 
with denominator dividing $\mu_F(n)^2$.
\item For every $n$ the number $p_n(a/F)$ is rational with
denominator dividing $\lcm(1,\ldots,n)^2\mu_F(n)^2$. 
\item For every $\epsilon>0$ we have that $|q_n(a/F)|,
|p_n(a/F)|<e^{\epsilon n}$ for sufficiently large $n$.
\item
Suppose $p^r||F$ where $r>0$ and $a$ is not divisible by $p$.
Then
$$|p_n(a/F)-\Theta_p(a/F)q_n(a/F)|_p\le p^2n^2p^{-2n(r+1/(p-1))}$$
for every $n$.
\end{enumerate}
\end{proposition}

{\bf Proof}.
The numbers $q_n(a/F)$ are given by
$$\sum_{k=0}^n {(1-a/F)_{n-k}(a/F)_k^2\over (n-k)!(k!)^2}.$$
The first assertion follows from Lemma \ref{pochhammerrational}.

The generating function of the $p_n(a/F)$ is the series $y_1$ with
$x=a/F$. To apply Proposition \ref{almostintegral} we replace
$z$ by $F^2\lambda^2 z$ and $x$ by $a/F$ in the equation $L_2(y)=0$
where $\lambda=\prod_{q|F}q^{1/(q-1)}$.
When we take the ring $R=\bbbz[q^{1/(q-1)}]_{q|F}$, the conditions of
Proposition \ref{almostintegral} are still satisfied with
$y_0(F^2\lambda^2 z)\in R[[z]]$ as power series solution.
From this Proposition it follows that the $n$-th coefficient of
$y_1(F^2\lambda^2 z)$ has denominator dividing $\lcm(1,\ldots,n)^2$.
Thus, our second statement follows.

The third statement on the Archimedean size of $q_n(a/F)$ and $p_n(a/F)$
follows from the fact that $y_0$ and $y_1$ have radius of convergence $1$.

The fourth statement follows from Proposition \ref{remainderI}.
It is a consequence of Lemma \ref{pochhammerrational} that 
$$\left|\fc{k}{a/F}\fc{k}{1-a/F}\right|_p<k^2p^{2-2k(r+1/(p-1))}$$
for all $k$. Hence 
$$|\Theta_p(n,a/F)|_p\le \max_{k\ge n}<k^2p^{2-2k(r+1/(p-1))}$$
from which our assertion follows. 

\qed
\medskip

We are now ready to state our irrationality results for $\Theta_p(a/F)$.

\begin{theorem}\label{mainI}
Let $a$ be an integer not divisible by $p$ and $F$ a natural number
divisible by $p$. Define $r$ by $|F|_p=p^{-r}$. Suppose that
$$\log F +\sum_{q|F}{\log q\over q-1} +1< 2r\log p +2{\log p\over p-1}.
\eqno{{\rm(A)}}$$
Then the $p$-adic number $\Theta_p(a/F)$ is irrational.
\end{theorem}

{\bf Proof}. Let $\epsilon>0$.
According to Proposition \ref{estimatesI}, $q_n(a/F),p_n(a/F)$
have a common
denominator dividing $Q_n:=\lcm(1,2,\ldots,n)^2\mu_F(n)^2$.
We also have, for $n$ large enough, $|q_n(a/F)|,|p_n(a/F)|<e^{n\epsilon}$. 
Furthermore $p_n(a/F)-q_n(a/F)\Theta_p(a/F)$ is non-zero for
infinitely many $n$. This follows from the fact that
$$p_{n+1}(x)q_n(x)-p_n(x)q_{n+1}(x)=1/(n+1)^2,$$
which can be shown by induction using recurrence (\ref{recurrenceI}).
We get
$$|Q_np_n(a/F)-Q_nq_n(a/F)\Theta_p(a/F)|_p
<p^{(-2r-2/(p-1)+\epsilon)n}|Q_n|_p$$
when $n$ is large enough.
We now apply Proposition \ref{irrationality} with $\alpha=\Theta_p(a/F),
q_n=Q_nq_n(a/F),p_n=Q_np_n(a/F)$. Notice that, for $n$ large enough,
\begin{eqnarray*}
|q_n|,|p_n|&<&e^{\epsilon n}\lcm(1,2,\ldots,n)^2\mu_F(n)^2\\
&<& \exp\left(n\epsilon+2n(1+\epsilon)+2n\log F+2n\sum_{q|F}{\log q\over q-1}
\right)
\end{eqnarray*}
In the latter we used the estimate $\lcm(1,\ldots,n)<e^{(1+\epsilon)n}$
which follows from the prime number theorem. Since 
$$|Q_n|_p<p^{-2rn}p^{-2[n/(p-1)]}\le p^{2-2n(r+1/(p-1))}$$
we get the estimate
$$|p_n-q_n\Theta_p(a/F)|_p<\exp\left(-4rn\log p-4n{\log p\over p-1}+
\epsilon n\right).$$
From Proposition \ref{irrationality} we can conclude irrationality
of $\Theta_p(a/F)$ if 
$$2+3\epsilon+2\log F+2\sum_{q|F}{\log q\over q-1}-4r\log p-4{\log p
\over p-1}+\epsilon<0.$$
From assumption (A) in our Theorem this certainly follows if $\epsilon$
is chosen sufficiently small. 

\qed
\medskip

\begin{corollary}\label{irrationalityI}
Let $\chi_8$ be the primitive even character modulo 8. Then
$\zeta_2(2),\zeta_3(2)$ and $L_2(2,\chi_8)$ are irrational.
\end{corollary}

{\bf Proof}. This is a direct consequence of Theorem \ref{mainI}
and Proposition \ref{cohen}. 

\qed
\medskip

Unfortunately condition (A) in Theorem \ref{mainI} is not good enough
to provide irrationality of $\Theta_3(1/6)$ which is related to $L_3(2,
\chi_{12})$. 

\section{Pad\'e approximations II}\label{padeII}
In \cite{Sti90} Stieltjes discovered the following continued fraction expansion
$$R(x)=\cfrac{-2}{2x-1+\cfrac{a_1}{2x-1+\cfrac{a_2}{\ddots}}}$$
with $a_n=n^4/(4n^2-1)$.

The convergents to this continued fraction were explicitly determined
by Touchard \cite{Tou56} and Carlitz \cite{Car57}. They were also used by Pr\'evost [Pre96] in his 
alternative irrationality proofs for $\zeta(2)$ and $\zeta(3)$.
In this section we give a self-contained derivation of the properties
of these convergents.

The numerator and denominator of the convergents satisfy the
recurrence relation
$$U_{n+1}=(2x-1)U_n+{n^4\over(4n^2-1)}U_{n-1}.$$
If we set $U_n=(n!)^2u_n/(1\cdot3\cdot5\cdots 2n-1)$ we get
\begin{equation}\label{recurrenceII}
(n+1)^2u_{n+1}=(2n+1)(2x-1)u_n+n^2u_{n-1}
\end{equation}

Consider the solutions $q_n(x)$ and $p_n(x)$ of this recurrence given
by
\begin{eqnarray*}
q_0(x) &=& 1\\
q_1(x) &=& 2x-1\\
q_2(x) &=& 3x^2-3x+1\\
q_3(x) &=& 10x^3/3-5x^2+11x/3-1\\
&&\cdots
\end{eqnarray*}

and
\begin{eqnarray*}
p_0(x) &=& 0\\
p_1(x) &=& -2\\
p_2(x) &=& -(6x-3)/2\\
p_3(x) &=& -(60x^2-60x+31)/18\\
&&\cdots
\end{eqnarray*}

Then $p_n(x)/q_n(x)$ are the convergents of our continued fraction.
To determine $q_n(x)$ we consider the generating function
$$y_0(z)=\sum_{n=0}^{\infty}q_n(x)z^n$$
and note that it satisfies the second order linear differential
equation
$$L_2(y)=(z^3-z)y''+(3z^2+(4x-2)z-1)y'+(z+2x-1)y$$
where the $'$ denotes differentation with respect to $z$.
At $z=0$ the equation $L_2(y)=0$ has a unique (up to a scalar factor)
holomorphic solution. In a straightforward manner one can thus
verify that
$$y_0(z)=(1+z)^{2x-1}\ _2F_1(x,x,1,z^2).$$
Comparison of coefficients of $z^n$ gives us the following explicit formula,
$$q_n(x)=\sum_{k\le n/2}{2x-1\choose n-2k}{-x\choose k}^2.$$

Consider the generating function of the $p_n(x)$,
$$y_1(z)=\sum p_n(x)z^n.$$
It is straightforward,
using the recurrence relation,
to see that $y_1$ satisfies the inhomogeneous equation $L_2(y_1)=2$.

We must now show that the rational functions $p_n(x)/q_n(x)$ approximate
$R(x)$ in $K$. To that end we define for each $n$,
$$R(n,x)=(-1)^n\sum_{k=0}^{\infty}{k(k-1)\cdots(k-n+1)\over
(k+1)(k+2)\cdots(k+n+1)}\fc{k}{x}.$$
Notice that it follows from Proposition \ref{basicidentities}
that $R(0,x)=-R(x)$.

\begin{proposition}\label{remainderII}
Letting notations be as above, we have for each $n$,
$$p_n(x)-q_n(x)R(x)=R(n,x)=O(1/x^{n+1})$$
as Laurent series in $1/x$.
\end{proposition}

{\bf Proof}. Letting
$$F(k,n)=(-1)^n{n!\over(k+1)\cdots(k+n+1)}{k\choose n}\fc{k}{x}$$
the Zeilberger algorithm shows that
$$-n^2F(k,n-1)-(2n+1)(2x-1)F(k,n)+(n+1)^2F(k,n+1)=
\Delta_k(2F(k,n)(x+k)(2n+1)).$$
When $n\ge1$ summation over $k$ yields
$$-n^2R(n-1,x)-(2n+1)(2x-1)R(n,x)+(n+1)^2R(n+1,x)=0.$$
When $n=0$ we get
$$-(2x-1)R(0,x)+R(1,x)=2.$$
From this, and the fact that $R(0,x)=-R(x)$ we conclude that
$$p_n(x)-q_n(x)R(x)=R(n,x)$$
for all $n\ge0$.

\qed
\medskip

There exist several interesting ways to write $q_n(x)$ as a
binomial sum, see for example \cite{Pre96}. One of them is
$$q_n(x)=(-1)^n\sum_{k=0}^n{n\choose k}{n+k\choose k}{-x\choose k}.$$
Taking $x=-n$ one recovers the Ap\'ery numbers for $\zeta(2)$ again.

Furthermore, in \cite{Pre96} we find the explicit expression
$$p_n(x)=(-1)^{n}\sum_{k=1}^n{n\choose k}{n+k\choose k}{-x\choose k}
\sum_{j=1}^k{(-1)^{j}\over j^2{-x\choose j}}.$$

\section{Application II}
In this section we prove irrationality for a large class of $p$-adic
numbers of the form $R_p(a/F)$ where $|a/F|_p>1$.

\begin{proposition}\label{estimatesII}
Let $p_n(x),q_n(x)$ be as in the previous section
and let $\mu_F(n)$ be as in Lemma 
\ref{pochhammerrational}.
Then, 
\begin{enumerate}
\item For every $n$ the number $q_n(a/F)$ is rational 
with denominator dividing $\mu_F(n)$.
\item For every $n$ the number $p_n(a/F)$ is rational with
denominator dividing $\lcm(1,\ldots,n)^2\mu_F(n)$. 
\item For every $\epsilon>0$ we have that $|q_n(a/F)|,
|p_n(a/F)|<e^{\epsilon n}$ for sufficiently large $n$.
\item
Suppose $p^r||F$ where $r>0$ and $a$ is not divisble by $p$,
we have
$$|p_n(a/F)-R_p(a/F)q_n(a/F)|_p\le (2n+1)\cdot p^{1-n(r+1/(p-1))}$$
for every $n$.
\end{enumerate}
\end{proposition}

{\bf Proof}.
The numbers $q_n(a/F)$ are given by
$$\sum_{k\le n/2} {2a/F-1\choose n-2k} {(a/F)_k^2\over (k!)^2}.$$
The first assertion follows from Lemma \ref{pochhammerrational}.

The generating function of the $p_n(a/F)$ is the series $y_1$ with
$x=a/F$. To apply Proposition \ref{almostintegral} we replace
$z$ by $F\lambda z$ and $x$ by $a/F$ in the equation $L_2(y)=0$
where $\lambda=\prod_{q|F}q^{1/(q-1)}$.
When we take the ring $R=\bbbz[q^{1/(q-1)}]_{q|F}$, the conditions of
Proposition \ref{almostintegral} are still satisfied with
$y_0(F\lambda z)\in R[[z]]$ as power series solution.
From this Proposition it follows that the $n$-th coefficient of
$y_1(F\lambda z)$ has denominator dividing $\lcm(1,\ldots,n)^2$.
Thus, our second statement follows.

The third statement on the Archimedean size of $q_n(a/F)$ and $p_n(a/F)$
follows from the fact that $y_0$ and $y_1$ have radius of convergence $1$.

The fourth statement follows from Proposition \ref{remainderII}.
It is a consequence of Lemma \ref{pochhammerrational} that 
$$\left|\fc{k}{a/F}\right|_p<kp^{1-k(r+1/(p-1))}$$
for all $k$. Moreover,
$$\left|{n!\over (k+1)(k+2)\cdots(k+n+1)}\right|_p=
\left|\sum_{l=0}^n
(-1)^l{n\choose l}{1\over k+l+1}\right|_p\le k+n+1.$$
Hence 
$$|R_p(n,a/F)|_p\le \max_{k\ge n}<(k+n+1)p^{1-k(r+1/(p-1))}$$
from which our assertion follows. 

\qed
\medskip

We are now ready to state our irrationality results for $R_p(a/F)$.

\begin{theorem}\label{mainII}
Let $a$ be an integer not divisible by $p$ and $F$ a natural number
divisible by $p$. Define $r$ by $|F|_p=p^{-r}$. Suppose that
$$\log F +\sum_{q|F}{\log q\over q-1} +2< 2r\log p +2{\log p\over p-1}.
\eqno{{\rm(B)}}$$
Then the $p$-adic number $R_p(a/F)$ is irrational.
\end{theorem}

{\bf Proof} Let $\epsilon>0$.
According to Proposition \ref{estimatesII}, $q_n(a/F),p_n(a/F)$
have a common
denominator dividing $Q_n:=\lcm(1,2,\ldots,n)^2\mu_F(n)$.
We also have, for $n$ large enough, $|q_n(a/F)|,|p_n(a/F)|<e^{n\epsilon}$. 
Furthermore $p_n(a/F)-q_n(a/F)R_p(a/F)$ is non-zero for
infinitely many $n$. This follows from the fact that
$$p_{n+1}(x)q_n(x)-p_n(x)q_{n+1}(x)=(-1)^{n-1}\cdot2/(n+1)^2.$$
This can be shown by induction using recurrence (\ref{recurrenceII}).
We get
$$|Q_np_n(a/F)-Q_nq_n(a/F)R_p(a/F)|_p
<p^{(-r-1/(p-1)+\epsilon)n}|Q_n|_p$$
when $n$ is large enough.
We now apply Proposition \ref{irrationality} with $\alpha=R_p(a/F),
q_n=Q_nq_n(a/F),p_n=Q_np_n(a/F)$. Notice that, for $n$ large enough,
\begin{eqnarray*}
|q_n|,|p_n|&<&e^{\epsilon n}\lcm(1,2,\ldots,n)^2\mu_F(n)\\
&<& \exp\left(n\epsilon+2n(1+\epsilon)+n\log F+n\sum_{q|F}{\log q\over q-1}
\right)
\end{eqnarray*}
In the latter we used the estimate $\lcm(1,\ldots,n)<e^{(1+\epsilon)n}$
which follows from the prime number theorem. Since 
$$|Q_n|_p<p^{-rn}p^{-[n/(p-1)]}\le p^{1-n(r+1/(p-1))}$$
we get the estimate
$$|p_n-q_n\Theta_p(a/F)|_p<\exp\left(-2rn\log p-2n{\log p\over p-1}+
\epsilon n\right).$$
From Proposition \ref{irrationality} we can conclude irrationality
of $R_p(a/F)$ if 
$$2+3\epsilon+\log F+\sum_{q|F}{\log q\over q-1}-2r\log p-2{\log p
\over p-1}+\epsilon<0.$$
From assumption (B) in our Theorem this certainly follows if $\epsilon$
is chosen sufficiently small. 

\qed
\medskip

\begin{corollary}
Let $p$ be a prime and $F$ a power of $p$ with $F\ne2$. Let $a$ be an
integer not divisible by $p$. Then $\omega(a)^{-1}H_p(2,a,F)$ is
irrational. For $F=2$ we have that $H_2(2,1,2)=0$.
\end{corollary}

{\bf Proof} Verify that condition (B) of Theorem \ref{mainII} holds for
every prime power $F>3$. The vanishing of $H_2(2,1,2)$ follows
from $R(x)+R(1-x)=0$ which implies $2R_2(1/2)=0$. Finally, irrationality
of $H_3(2,1,3)=H_3(2,2,3)$ follows from the irrationality of $\zeta_3(2)$
proved in Corollary \ref{irrationalityI}.

\qed
\medskip

\section{Pad\'e approximations III}\label{padeIII}
In this section we prove the following continued fraction
expansion of $T(x)=\sum_{n=0}^{\infty}(n+1)B_n(-1/x)^{n+2}$,
namely
$$T(x)=\cfrac{1}{a_1-\cfrac{1^6}{a_2 -\cfrac{2^6}{a_3 -\ddots}}}$$
where $a_n=(2n-1)(2x^2-2x+n^2-n+1)$. 
In \cite{Sti90}{(23)]}
we find a related continued fraction for $x^2T(x)-1-1/x$, but
we prefer the one presented here because it has simpler properties.
Moreover, by substituting $x=1/4$, we obtain a continued fraction
expansion which converges rapidly $2$-adically to $T_2(1/4)=64\zeta_2(3)$.
Without proof we note that Calegari's approximations 
$a_n/b_n$ to $\zeta_2(3)$ (see proof of \cite[Thm 3.4]{Cal05} 
coincide with the fractions $-{p_n(1/4)-p_{n-1}(1/4)
\over 64(q_n(1/4)-q_{n-1}(1/4))}$. Here $p_n(x)/q_n(x)$ are the
convergents of the continued fraction for $T(x)$, to be specified
below. A similar remark holds for Calegari's approximations to 
$\zeta_3(3)$.

Our study of the convergents of the continued fraction expansion begins
with the observation that the numerators and demoninators of the 
convergents can be normalised in such a way that they are solutions 
of the recurrence
$$U_{n+1}=(2n+1)(2x^2-2x+n^2+n+1)U_n-n^6U_{n-1}.$$
Replace $U_n$ by $n^3u_n$ to find
\begin{equation}\label{recurrenceIII}
(n+1)^3u_{n+1}=(2n+1)(2x^2-2x+n^2+n+1)u_n-n^3u_{n-1}
\end{equation}
Two independent solutions $q_n(x),p_n(x)$ are given by
\begin{eqnarray*}
q_0(x)&=&1\\
q_1(x)&=&2x^2-2x+1\\
q_2(x)&=&(3x^4-6x^3+9x^2-6x+2)/2\\
q_3(x)&=&(10x^6-35x^5+85x^4-120x^3+121x^2-66x+18)/18\\
&&\cdots
\end{eqnarray*}
and
\begin{eqnarray*}
p_0(x)&=&0\\
p_1(x)&=&2\\
p_2(x)&=&3(2x^2-2x+3)/4\\
p_3(x)&=&(60x^4-120x^3+360x^2-300x+251)/108\\
&&\cdots
\end{eqnarray*}

Consider the generating function 
$$Y_0(z)=\sum_{n=0}^{\infty}q_n(x)z^n.$$
Using the recurrence we see that $Y_0(z)$ is solution
of the linear differential equation
\begin{eqnarray*}
L_3(y)&=&z^2(z-1)^2y'''+3z(2z-1)(z-1)y''+\\
&&+(7z^2-(4x^2-4x+8)z+1)y'+(z-2x^2+2x-1)y
\end{eqnarray*}

One can verify in a straightforward 
manner that this equation is the symmetric square of the second order equation
$$L_2(y)=z(z-1)^2y''+(z-1)(2z-1)y'+(z/4+x-x^2-1/2)y.$$

The unique power series solution in $z$ of $L_2(y)=0$ reads
$$y_0(z)=(1-z)^{x-1/2}\ _2F_1(x,x,1,z).$$
As a consequence the function $Y_0(z)$ equals
$$Y_0(z)=y_0(z)^2=(1-z)^{2x-1}\ _2F_1(x,x,1,z)^2.$$
By comparison of coefficients we would be able to compute an
explicit expression for $q_n(x)$. But this would be a double
summation. 
A much nicer expression for $q_n(x)$ can be found from \cite{Pre96}. It reads
$$q_n(x)=\sum_{k=0}^n{n \choose k}{n+k\choose k}{-x\choose k}
{-x+k\choose k}.$$
Notice that $x=-n$ recovers the Ap\'ery numbers for $\zeta(3)$.

Let $Y_1$ be the generating function of the $p_n(x)$. Then $Y_1$
satisfies the inhomogeneous equation $L_3(Y_1)=1$.

An explicit formula from \cite{Pre96} reads
$$p_n(x)=\sum_{k=1}^n{n\choose k}{n+k\choose k}
\sum_{j=1}^k{(-1)^{j-1}\over j^3}\ {{k-x\choose k-j}{-x-j\choose k-j}\over
{k\choose j}^2}.$$

We must now show that the rational functions $p_n(x)/q_n(x)$ approximate
$T(x)$ in $K$. To that end we define for each $n$,
$$T(n,x)=(-1)^n\sum_{k=0}^{\infty}{k(k-1)\cdots(k-n+1)\over
(k+1)(k+2)\cdots(k+n+1)}\fc{k}{x}\fc{k}{1-x}.$$
Notice, using Proposition \ref{basicidentities}, that $T(0,x)=-T(x)$.

\begin{proposition}\label{remainderIII}
Letting notations be as above, we have for each $n$,
$$p_n(x)-q_n(x)T(x)=T(n,x)=O(1/x^{2n+2})$$
as Laurent series in $1/x$.
\end{proposition}

{\bf Proof}. Letting
$$F(k,n)=(-1)^n{n!\over(k+1)\cdots(k+n+1)}{k\choose n}\fc{k}{x}
\fc{k}{1-x}$$
the Zeilberger algorithm shows that
$$n^3F(k,n-1)-(2n+1)(2x^2-2x+n^2+n+1)F(k,n)+(n+1)^3F(k,n+1)$$
$$=
\Delta_k(2F(k,n)(x+k)(1-x+k)(2n+1)).$$
When $n\ge1$ summation over $k$ yields
$$n^3T(n-1,x)-(2n+1)(2x^2-2x+n^2+n+1)T(n,x)+(n+1)^3T(n+1,x)=0.$$
When $n=0$ we get
$$-(2x^2-2x+1)T(0,x)+T(1,x)=2.$$
From this, and the fact that $T(0,x)=-T(x)$ we conclude that
$$p_n(x)-q_n(x)T(x)=T(n,x)$$
for all $n\ge0$.

\qed
\medskip

\section{Application III}
In this section we prove irrationality for a large class of $p$-adic
numbers of the form $T_p(a/F)$ where $|a/F|_p>1$.

\begin{proposition}\label{estimatesIII}
Let notations be as above and let $\mu_F(n)$ be as in Lemma 
\ref{pochhammerrational}.
Then, 
\begin{enumerate}
\item For every $n$ the number $q_n(a/F)$ is rational 
with denominator dividing $\mu_F(n)^2$.
\item For every $n$ the number $p_n(a/F)$ is rational with
denominator dividing $\lcm(1,\ldots,n)^3\mu_F(n)^2$. 
\item For every $\epsilon>0$ we have that $|q_n(a/F)|,
|p_n(a/F)|<e^{\epsilon n}$ for sufficiently large $n$.
\item
Suppose $p^r||F$ where $r>0$ and $a$ is not divisible by $p$,
we have
$$|p_n(a/F)-T_p(a/F)q_n(a/F)|_p\le (2n+1)n^2\cdot p^{2-2n(r+1/(p-1))}$$
for every $n$.
\end{enumerate}
\end{proposition}

{\bf Proof}.
The numbers $q_n(a/F)$ are the coefficients of $(1-z)^{2a/F-1}
\ _2F_1(a/F,a/F,1,z)^2$. Let again, $\lambda=\prod_{q|F}q^{1/(q-1)}$
and let $R$ be the ring of integers in $\bbbq(\lambda)$. 
Then, by Lemma \ref{pochhammerrational}, we have 
$(1-\lambda^2z)^{2a/F-1},\ _2F_1(a/F,a/F,1,\lambda^2z)\in R[[z]]$.
Hence part i) follows.

The generating function of the $p_n(a/F)$ is the series $Y_1$ with
$x=a/F$. To apply Proposition \ref{almostintegral} we replace
$z$ by $F\lambda z$ and $x$ by $a/F$ in the equation $L_2(y)=0$.
The conditions of
Proposition \ref{almostintegral} are still satisfied with
$y_0(F\lambda^2 z)\in R[[z]]$ as power series solution.
From this Proposition it follows that the $n$-th coefficient of
$Y_1(F\lambda^2 z)$ has denominator dividing $\lcm(1,\ldots,n)^3$.
Thus, our second statement follows.

The third statement on the Archimedean size of $q_n(a/F)$ and $p_n(a/F)$
follows from the fact that $y_0$ and $y_1$ have radius of convergence $1$.

The fourth statement follows from Proposition \ref{remainderIII}.
It is a consequence of Lemma \ref{pochhammerrational} that 
$$\left|\fc{k}{a/F}\fc{k}{1-a/F}\right|_p<k^2p^{2-2k(r+1/(p-1))}$$
for all $k$. Moreover,
$$\left|{n!\over (k+1)(k+2)\cdots(k+n+1)}\right|_p=
\left|\sum_{l=0}^n
(-1)^l{n\choose l}{1\over k+l+1}\right|_p\le k+n+1.$$
Hence 
$$|T_p(n,a/F)|_p\le \max_{k\ge n}<(k+n+1)k^2 p^{2-2k(r+1/(p-1))}$$
from which our assertion follows. 

\qed
\medskip

We are now ready to state our irrationality results for $T_p(a/F)$.

\begin{theorem}\label{mainIII}
Let $a$ be an integer not divisible by $p$ and $F$ a natural number
divisible by $p$. Define $r$ by $|F|_p=p^{-r}$. Suppose that
$$\log F +\sum_{q|F}{\log q\over q-1} +3/2< 2r\log p +2{\log p\over p-1}.
\eqno{{\rm(C)}}$$
Then the $p$-adic number $T_p(a/F)$ is irrational.
\end{theorem}

{\bf Proof} Let $\epsilon>0$.
According to Proposition \ref{estimatesIII}, the rational
numbers $q_n(a/F),p_n(a/F)$
have a common
denominator dividing $Q_n:=\lcm(1,2,\ldots,n)^3\mu_F(n)^2$.
We also have, for $n$ large enough, $|q_n(a/F)|,|p_n(a/F)|<e^{n\epsilon}$. 
Furthermore $p_n(a/F)-q_n(a/F)T_p(a/F)$ is non-zero for
infinitely many $n$. This follows from the fact that
$$p_{n+1}(x)q_n(x)-p_n(x)q_{n+1}(x)=1/(n+1)^3.$$
This can be shown by induction using recurrence (\ref{recurrenceII}).
We get
$$|Q_np_n(a/F)-Q_nq_n(a/F)T_p(a/F)|_p
<p^{2n(-r-1/(p-1)+\epsilon)}|Q_n|_p$$
when $n$ is large enough.
We now apply Proposition \ref{irrationality} with $\alpha=T_p(a/F),
q_n=Q_nq_n(a/F),p_n=Q_np_n(a/F)$. Notice that, for $n$ large enough,
\begin{eqnarray*}
|q_n|,|p_n|&<&e^{\epsilon n}\lcm(1,2,\ldots,n)^3\mu_F(n)^2\\
&<& \exp\left(n\epsilon+3n(1+\epsilon)+2n\log F+2n\sum_{q|F}{\log q\over q-1}
\right)
\end{eqnarray*}
In the latter we used the estimate $\lcm(1,\ldots,n)<e^{(1+\epsilon)n}$
which follows from the prime number theorem. Since 
$$|Q_n|_p<p^{-2rn}p^{-2[n/(p-1)]}\le p^{2-2n(r+1/(p-1))}$$
we get the estimate
$$|p_n-q_nT_p(a/F)|_p<\exp\left(-4rn\log p-4n{\log p\over p-1}+
\epsilon n\right).$$
From Proposition \ref{irrationality} we can conclude irrationality
of $T_p(a/F)$ if 
$$3+4\epsilon+2\log F+2\sum_{q|F}{\log q\over q-1}-4r\log p-4{\log p
\over p-1}+\epsilon<0.$$
From assumption (C) in our Theorem this certainly follows if $\epsilon$
is chosen sufficiently small. 

\qed
\medskip

\begin{corollary}
Let $p$ be a prime and $F$ a power of $p$ with $F>2$. Let $a$ be an
integer not divisible by $p$. Then $\omega(a)^{-2}H_p(3,a,F)$ is
irrational. 
\end{corollary}

{\bf Proof} Verify that condition (C) of Theorem \ref{mainIII} holds for
every prime power $F>2$.

\qed
\medskip

\begin{corollary}
Let $\chi_d$ be a primitive even character modulo $d$.
Then the following numbers are irrational: $\zeta_2(3),\zeta_3(3),
\zeta_5(3)-L_3(3,\chi_5),\zeta_5(3)+L(3,\chi_5),\zeta_2(3)-
L_2(3,\chi_8),\zeta_2(3)+L_2(3,\chi_8)$.
\end{corollary}

{\bf Proof} Use the previous Corollary and Proposition
\ref{valuesat3}.

 \qed
\medskip

\section{Pad\'e approximations IV}
So far we have studied continued fraction expansions of the 
functions $\Theta(x),R(x)$ and $T(x)$. Clearly $R(x),T(x)$ are the
generator series of the Bermoulli numbers and its derivatives.
We like to remark here that the coefficients $(2^{n+1}-2)B_n$ of
$\Theta(x)$ are actually $(n-1)T_{n-1}$ for $n>1$ where $T_n$
is the hyperbolic tangent number defined by
$$\tanh(t/2)={e^t-1\over e^t+1}=\sum_{n=0}^{\infty}{T_n\over n!}t^n.$$
This follows from the observation that 
$$\sum_{n=0}^{\infty}(2^{n+1}-2){B_n\over n!}t^n={4t\over e^{2t}-1}-
{2t\over e^t-1}=t\tanh(t/2)-t.$$
The series $\Theta(x)$ is also related to the Euler numbers $E_n$
via
$$\Theta((1-z)/2)=-4\sum_{n=0}^{\infty}(n+1)E_n(1/z)^{n+2}.$$
The Euler numbers are defined by

The only interesting additional continued fraction ($S$-fraction in the
sense of Stieltjes) we have been able to find is one for
$$\theta(x)=\sum_{n=1}^{\infty}(2^{n+1}-2){B_n\over n}(-1/x)^n.$$
It reads
$$\theta(x)=\cfrac{2}{2x-1+\cfrac{1}{2x-1+\cfrac{4}{2x-1+\cfrac{9}{2x-1+
\cfrac{16}{2x-1+\ddots}}}}}.$$
We will not give any proofs here (they are parallel to the previous
sections), but only quote some formulas. The recurrence relation
involved with this continued fraction is
$$U_{n+1}=(2x-1)U_n+n^2U_{n-1}$$
Substitute $U_n=n!u_n$. Then,
$$(n+1)u_{n+1}=(2x-1)u_n+nu_{n-1}.$$
Consider the solutions
\begin{eqnarray*}
q_0(x)&=&1\\
q_1(x)&=&2x-1\\
q_2(x)&=&2x^2-2x+1\\
q_3(x)&=&(2x-1)(2x^2-2x+3)/3\\
&&\cdots
\end{eqnarray*}
and
\begin{eqnarray*}
p_0(x)&=&0\\
p_1(x)&=&2\\
p_2(x)&=&2x-1\\
p_3(x)&=&(4x^2-4x+5)/3\\
&&\cdots
\end{eqnarray*}
The generating function 
$$y_0(z)=\sum_{n=0}^{\infty}q_n(x)z^n$$
satisfies the differential equation
$$(1-z^2)y'-(2x-1+z)y=0.$$
One easily recovers that
$$y_0(z)=(1-z)^{-x}(1+z)^{x-1}.$$
From this we infer with a bit of effort
$$q_n(x)=(-1)^n\sum_{k=0}^n{n\choose k}{-x\choose k}2^k.$$
The generating function
$$y_1(z)=\sum_{n=0}^{\infty}p_n(x)z^n$$
satisfies
$$(1-z^2)y'-(2x-1+z)y=2.$$

We also have the identity
$$\theta(x)=\sum_{k=0}^{\infty}{1\over 2^k}\fc{k}{x}.$$
Let us define 
$$\theta(n,x)=(-1)^n\sum_{k=0}^{\infty}{{k\choose n}\over 2^k}\fc{k}{x}.$$
Then we have the Pad\'e approximation property
$$p_n(x)-q_n(x)\theta(x)=-\theta(n,x)=O(1/x^n).$$
Just as in the previous sections we could apply this to $p$-adic irrationality
proofs, but will not pursue this here. We only remark that $L_2(1,\chi_8)$
and $L_3(1,\chi_{12})$ can be proven irrational.

\section{Appendix: $p$-Adic Hurwitz series}

Let $p$ be a prime. Let 
$F$ be a positive integer and $a$ an integer not divisible by $F$.
Later, when we define $p$-adic functions, we shall assume 
in addition that $p$
divides $F$.
Define the Hurwitz zeta-function 
$$H(s,a,F)=\sum_{n=0}^{\infty}{1\over (a+nF)^s}.$$
This series converges for all $s\in\bbbc$ with real part $>1$.
As is well-known $H(s,a,F)$ can be continued analytically to the
entire complex $s$-plane, with the exception of a pole at $s=1$.

Let $n$ be an integer $\ge1$. To determine the value $H(1-n,a,F)$
we expand
$${te^{at}\over e^{Ft}-1}=\sum_{n\ge0}{B_{n}(a,F)
\over n!}t^n.$$
Then we have

\begin{proposition}\label{specialvalue}
For any $F,a,n$ we have
$$H(1-n,a,F)=-B_n(a,F)/n.$$
\end{proposition}

We can express $B_n(a,F)$ in terms of the ordinary Bernoulli numbers
$B_k$ which are given by
$${t\over e^t-1}=\sum_{k\ge0}{B_k\over k!}t^k.$$
We get

\begin{lemma}
For any positive integer $n$,
$$B_{n}(a,F)={a^n\over F}\sum_{j=0}^n{n\choose j}B_j
\left({F\over a}\right)^j.$$
\end{lemma}

{\bf Proof}. We expand in powers of $t$,
\begin{eqnarray*}
{te^{at}\over e^{Ft}-1}&=&
{1\over F}e^{at}\sum_{j\ge0}{B_j\over j!}(Ft)^j\\
&=&{1\over F}\sum_{n\ge0}\sum_{i+j=n}
\left({a^i\over i!}{B_j\over j!}F^j\right)t^n\\
&=&{1\over F}\sum_{n\ge0}\left(\sum_{j=0}^n{n\choose j}
B_j\left({F\over a}\right)^j\right){a^n\over n!}t^n
\end{eqnarray*}
The proof of our Lemma now follows by comparison of the coefficient of
$t^n$.

\qed
\medskip

From now on we assume that $F$ is divisible by $p$ and $a$ is
not divisible by $p$. Then
$|B_j(F/a)^j|_p\to0$ as $j\to\infty$ and we can think of $p$-adic
interpolation.
We would like to interpolate the values $H(1-n,a,F)$ $p$-adically in $n$.
Strictly 
speaking this is impossible, but we can interpolate 
$H(1-n,a,F)\omega(a)^{-n}$
where $\omega$ is the Teichm\"uller character $\bbbz\to\bbbz_p$ given as
follows. When $p|m$ we define $\omega(m)=0$. When $\gcd(p,m)=1$ and
$p$ is odd, we define $\omega(m)^{p-1}=1$ and $\omega(m)\is m\mod{p}$.
When $p=2$ and $m$ odd, we define $\omega(m)=(-1)^{(m-1)/2}$.
We also define $<x>=\omega(x)^{-1}x$ for all integers $x$ not divisible
by $p$. Notice that $s\mapsto <x>^s$ is $p$-adically analytic on $\bbbz_p$.

We define the $p$-adic function $H_p$ by
$$H_p(s,a,F)={1\over F(s-1)}<a>^{1-s}\sum_{j=0}^{\infty}
{1-s\choose j}B_j\left({F\over a}\right)^j$$
for all $s\in\bbbz_p$.

Note in particular the value at $s=2$. This equals
$$H_p(2,a,F)=-{\omega(a)\over F^2} \sum_{j=0}^{\infty}B_j
\left(-{F\over a}\right)^{j+1}.$$
The latter summation is precisely the series $R(x)$, defined in the text
in which we have substituted $x=a/F$.

Finally we define the $p$-adic Kubota-Leopoldt $L$-series.
Let $\phi:\bbbz\to\bbba$ be a periodic function with period $f$.
Let $F=\lcm(f,p)$ if $p$ is odd and $F=\lcm(f,4)$ if $p=2$.
We now define the $p$-adic
$L$-series associated to $\phi$ by
$$L_p(s,\phi)=\sum_{a=1,a\not\equiv0\mod{p}}^F \phi(a)H_p(s,a,F).$$
We remark that the value of $L_p(s,\phi)$ remains the same
if we choose instead of $F$ a multiple period $mF$.
To see this it suffices to show that for all integers $n\ge0$,
$$\sum_{a=1,a\not\equiv0\mod{p}}^F\phi(a)H(1-n,a,F)=
\sum_{a=1,a\not\equiv0\mod{p}}^{mF}\phi(a)H(1-n,a,mF).$$
This follows from Proposition \ref{specialvalue} and the identity
$$\sum_{a=1,a\not\equiv0\mod{p}}^{mF}{t\phi(a)e^{at}
\over e^{mF}-1}=
\sum_{a=1,a\not\equiv0\mod{p}}^F{t\phi(a)e^{at}
\over e^{F}-1}$$
The latter follows from the periodicity
of $\phi$ with period $F$ and $p|F$.

When $\phi(n)=1$ for all $n$ we get the $p$-adic zeta-function
$$\zeta_p(s)=\sum_{a=1}^{p-1}H_p(s,a,p)$$
when $p$ is odd and when $p=2$,
$$\zeta_2(s)=H_2(s,1,4)+H_2(s,3,4).$$

\end{document}